\documentclass[12pt,a4paper]{article}
\usepackage{amsmath}
\usepackage{amssymb}
\usepackage{latexsym}
\usepackage{amsmath}
\usepackage{graphicx}
\usepackage{amsfonts}
\usepackage{amssymb}
\usepackage{amscd}

\begin{document}
\newtheorem{theorem}{Theorem}[section]
\newtheorem{remark}[theorem]{Remark}
\newtheorem{mtheorem}[theorem]{Main Theorem}
\newtheorem{bbtheo}[theorem]{The Strong Black Box}
\newtheorem{observation}[theorem]{Observation}
\newtheorem{proposition}[theorem]{Proposition}
\newtheorem{lemma}[theorem]{Lemma}
\newtheorem{testlemma}[theorem]{Test Lemma}
\newtheorem{mlemma}[theorem]{Main Lemma}
\newtheorem{note}[theorem]{}
\newtheorem{steplemma}[theorem]{Step Lemma}
\newtheorem{corollary}[theorem]{Corollary}
\newtheorem{notation}[theorem]{Notation}
\newtheorem{example}[theorem]{Example}
\newtheorem{definition}[theorem]{Definition}
\renewcommand{\labelenumi}{(\roman{enumi})}
\newcommand{\dach}[1]{\hat{\vphantom{#1}}}
\newcommand{\cp}{\widehat}
\newcommand{\dsum}{\bigoplus}
\newcommand{\pure}{\subseteq_\ast}
\numberwithin{equation}{section}
\def\Z{{ \mathbb Z}}
\def\E{{ \mathbb E}}
\def\H{{ \mathbb H}}
\def\N{{ \mathbb N}}
\def\S{{ \mathbb S}}
\def\DD{ \mathbb D}
\def\GG{ \mathbb G}
\def\Proof{{\sl Proof.}\quad}
\newcommand{\supp}[1]{\mbox{$\left[ \, #1\, \right]$}}
\newcommand{\suppl}[1]{\mbox{$\left[ \, #1\, \right]_\lambda$}}
\newcommand{\suppa}[1]{\mbox{$\left[ \, #1\, \right]_A$}}
\newcommand{\norm}[1]{\mbox{$\parallel\! #1 \!\parallel$}}
\newcommand{\norma}[1]{\mbox{$\parallel\! #1 \!\parallel_A$}}

\newcommand{\union}{\bigcup}
\newcommand{\card}[1]{\mbox{$\left| #1 \right|$}}
\def\black{\ {\hbox{\vrule width 4pt height 4pt depth
0pt}}}
\def\fine{\ \black\vskip.4truecm}
\def\BBB{ B_\BB}
\def\D{{\hat{D}}}
\def\Q{{\mathbb Q}}
\def\bS{{\mathbb S}}
\def\G{\hat{G}}
\def\T{{\cal T}}
\def\FF{\mathfrak  F}
\def\FFbar{\overline{\FF}}
\def\FFs{\mathfrak  F'}
\def\KK{\mathfrak  K}
\def\AA{\mathfrak  A}
\def\FFF{\langle \FF \rangle}
\def\FFFs{\langle \FFs \rangle}
\def\PP{{\mathfrak  P}}
\def\pp{{\mathfrak  p}}
\def\x{{\mathfrak  x}}
\def\y{{\mathfrak  y}}
\def\z{{\mathfrak  z}}
\def\C{{\mathfrak C}}
\def\X{{\mathfrak X}}
\def\Y{{\mathfrak Y}}
\def\RX{R\langle X \rangle}
\def\GX{\langle X \rangle}
\def\KKf{{$\mathfrak K_f$}}
\def\R{\widehat{R}}
\def\wZ{\widehat{\Z}}
\def\F{\widehat{F}}
\def\G{\widehat{G}}
\def\T{{\cal T}}
\def\B{\widehat{B}}
\def\restr{\restriction}
\def\cf{{\rm cf\,}}
\def\df{{\rm df\,}}
\def\id{{\rm id\,}}
\def\otp{{\rm otp\,}}
\def\br{{\rm br\,}}
\def\inf{{\rm inf\,}}
\def\sup{{\rm sup\,}}
\def\Br{{\rm Br\,}}
\def\Yphi{Y_{[\phi]}}
\def\Ypsi{Y_{[\psi]}}
\def\Xphi{X_{\tilde{\phi}}}
\def\Xpsi{X_{\tilde{\psi}}}
\def\a{\alpha}
\def\o{\omega}
\def\e{\varepsilon}
\def\hv{\widehat\varphi}
\def\v{\varphi}
\def\s{\sigma}
\def\k{\kappa}
\def\ale{\aleph_1}
\def\aln{\aleph_0}
\def\l{\lambda}
\def\b{\beta}
\def\g{\gamma}
\def\abar{\overline{\alpha}}
\def\vbar{\overline{\v}}
\def\Gbar{\overline{G}}
\def\Xbar{\overline{X}}
\def\xbar{\overline{x}}
\def\thetabar{\overline{\theta}}

\def\gbar{\overline{g}}
\def\aa{{\bf a}}
\def\to{\rightarrow}
\def\arr{\longrightarrow}
\def\restr{\upharpoonright}
\def\End{{\rm End\,}}
\def\Mon{{\rm Mon\,}}
\def\Mont{{\rm Mon_t\,}}
\def\pmont{{\rm pMon_t\,}}
\def\Ult{{\rm Ult\,}}
\def\Ines{{\rm Ines\,}}
\def\Hom{{\rm Hom\,}}
\def\Fin{{\rm Fin\,}}
\def\Ext{{\rm Ext}\,}
\def\Hom{{\rm Hom}\,}
\def\End{{\rm End}\,}
\def\pAut{{\rm pAut\,}}
\def\Aut{{\rm Aut\,}}
\def\Map{{\rm Map\,}}
\def\Im{{\rm Im\,}}
\def\ker{{\rm ker}\,}
\def\Ann{{\rm Ann}\,}
\def\defe{{\rm def}\,}
\def\rk{{\rm rk}\,}
\def\crk{{\rm crk}\,}
\def\nuc{{\rm nuc}\,}
\def\Dom{{\rm Dom}\,}
\def\Im{{\rm Im}\,}
\def\Yphi{Y_{[\phi]}}
\def\Ypsi{Y_{[\psi]}}
\def\Xphi{X_{\tilde{\phi}}}
\def\Xpsi{X_{\tilde{\psi}}}
\def\abar{\overline{\alpha}}
\def\ra{\rightarrow}
\def\arr{\longrightarrow}
\def\iff{\Longleftrightarrow}
\def\mm{{\mathfrak m}}
\def\X{{\mathfrak X}}
\def\Diam{\diamondsuit}
\title{{\sc Uniquely Transitive Torsion-free Abelian Groups}
\footnotetext{This work is supported by the project No.
I-706-54.6/2001 of the German-Israeli
Foundation for Scientific Research \& Development\\
AMS subject classification:\\
primary: 13C05, 13C10, 13C13, 20K15,20K25,20K30  \\
secondary: 03E05, 03E35\\
Key words and phrases: automorphism groups of torsion-free
abelian groups, \\
\\ GbSh 650 in Shelah's list of publications
}}

\author{R\"udiger G\"obel and Saharon Shelah}
\date{}
\maketitle

\begin{abstract}  We will answer a question raised by Emmanuel
Dror Farjoun concerning the existence of torsion-free abelian
groups $G$ such that for any ordered pair of pure elements there
is a unique automorphism mapping the first element onto the second
one. We will show the existence of such a group of cardinality
$\l$ for any successor cardinal $\l =\mu^+$ with $\mu =
\mu^{\aln}$.
\end{abstract}

\section{Introduction}

We will consider the set $ \pp G$ of all non-zero pure elements of
a torsion-free abelian group $G$. Recall that $g\in G$ is pure if
the equations $xn = g$ for natural numbers $n \neq 1$ have no
solution $x\in G$. Clearly every element of the automorphism group
$\Aut G$ of $G$ induces a permutation on the set $\pp G$ and it is
natural to consider groups where the action of $\Aut G$ on $\pp G$
is transitive: for any pair $x,y\in\pp G$ there is an automorphism
$\v\in\Aut G$ such that $x\v = y$.  In this case we will say that
$G$ is transitive, for short $G$ is a T-group. (Transitive groups
are $A$-transitive groups in Dugas, Shelah \cite{DS}). This kind
of consideration is well-known for abelian $p$-groups. It was
stimulated by Kaplansky and studied in many papers, see
\cite{Ka,Mg,Co} for instance. If $G$ is a free abelian group with
two pure elements $x$ and $y$, then there are two sets $X$ and $Y$
of free generators of $G$ such that $x\in X$ and $y\in Y$. We can
chose a bijection $\a: X\to Y$ with $x\a = y$ which extends to an
automorphism $\a \in \Aut G$. Hence free groups are T-groups and a
similar argument holds for a wide class of abelian groups. There
are T-groups like $\Z^{\ale}/\Z^{[\ale]}$ with $\Z^{[\ale]}$ the
set of all elements in $\Z^{\ale}$ of countable support, which are
$\ale$-free but not separable, see Dugas, Hausen \cite{DH}. The
existence of $\ale$-free, indecomposable T-groups in L for any
regular, not weakly compact cardinality was also shown in Dugas,
Shelah \cite [Theorem (b), p. 192] {DS}. This was used to answer a
problem in Hausen \cite{Ha1}, see also \cite{Ha2} and
\cite[Theorem (a), p. 192]{DS}.

Thus we may strengthening the action of $\Aut G$ on $\pp G$ and
say that $G$ is a U-group if any two automorphisms $\v, \v'\in
\Aut G$ with $g\v = g\v'$ for some $g\in \pp G$ must be the same
$\v =\v'$. Hence $G$ is a UT-group if and only if $G$ is both a
T-group as well as a U-group. Thus $G$ is a UT-group if and only
if $\Aut G$ is transitive and (every non-trivial automorphism
acts) fix point-free on $\pp G$. In connection with permutation
groups such action is also called `sharply transitive'. Note that
$\pp G$ may be empty, if $G$ is divisible for instance. In order
to avoid trivial cases we also require that $0\ne G\ne \Z$ and $G$
is of type $0$, hence $G$ is torsion-free and every element of $G$
is a multiple of an element in $\pp G$. If $G$ is of type $0$ and
not finitely generated then $|\pp G| =|G|$ is large and the
problem about the existence of UT-groups becomes really
interesting. This question is related to problems in homotopy
theory and was raised by Dror Farjoun. In response we want to show
the following
\begin{theorem} \label{super} For any successor cardinal
$\l = \mu^+$ with $\mu = \mu^{\aln}$ there is an $\ale$-free
abelian UT-group of cardinality $\l$.
\end{theorem}
We will also determine the endomorphism rings of these groups.
They are isomorphic to integral group rings $R= \Z F$ of groups
$F$ freely generated (as a non-abelian group) by $\l$ elements
(with $\l$ as in the theorem). Since endomorphisms of a group $G$
will act on the right (in accordance with used results from
\cite{GT,GW}),we will also view $G$ as a right $R$-module (and as
a left or right $\Z$-module). Using classical results on group
rings it will follow that $\Aut G = \pm F$, where $-1\in \Z$ is
scalar multiplication by $-1$, hence $\pm F$ is a direct product
of a group of order $2$ and $F$. Moreover $\Z F$ has no
idempotents except $0$ and $1$, hence $G$ in the theorem is
indecomposable and obviously the center of $\Z F$ is just $\Z$,
hence $\Z$ is also the center of $\End G$. Therefore Theorem
\ref{super} strengthens the Theorem in Dugas, Shelah \cite[p.
192]{DS} substituting T-groups by UT-groups and removing the
restriction $V=L$ to the constructible universe. Also note that it
is straightforward to replace the ground ring $\Z$ by a
$p$-reduced domain $S$ for some prime $p$.  Hausen's \cite{Ha2}
problem can be answered also in ordinary set theory. Further
applications can be found in Section \ref{Discus}.

First we would like to explain why constructing UT-groups is a
hard task, much harder then finding suitable T-groups. Because
$R^+$ above is a free abelian group, we can easily find groups $G$
with $R =\End G$, see \cite{DG,CG}. But there are still two
obstacles which must be taken into consideration. Often $|R| <
|G|$ in realization theorems, thus the units of $R$ which
represent $\Aut G$ will never act transitive on a bigger group and
$G$ can't be a T-group. More importantly, inspecting the
constructions in \cite{DG,CG}, it is clear that they provide no
control about the action, thus both U and T for UT-groups are a
problem. Note that in many earlier constructions $G$ has a free
dense and pure $R$-submodule of rank $>1$ mostly of rank $|G|$.
But T-groups must obviously be cyclic over their endomorphism
ring $R$, hence \cite{DG,CG} do not apply in principle. Inspecting
the proof in \cite{DS} it is easy to see that $G$ is not
torsion-free over its endomorphism ring. This comes from the list
of new variables $x,y,...$ added to $R$ in the construction in
order to make $R$ acting transitive on all pairs of pure
elements. Even refining the list of pure pairs in \cite{DS} it
seems hard to avoid clashes of related pairs such that $x-y$ for
example has a proper annihilator. Thus the groups in \cite{DS}
are T-groups and not UT-groups (even modifying the arguments).

Thus a new approach is need, which will be established in Section
3. We will use a geometric argument choosing carefully new partial
automorphisms for making $G$ transitive but with very small domain
and image in order to preserve the U-property for the new monoid.
Then we will feed the partial maps with pushouts to grow them up
and become real automorphism without destroying the UT-property.
At the end we will have a suitable subgroup $F$ of automorphisms
of some group $G$, thus $G$ becomes an $R$-module over the ring
$\Z F =^\df R$.

Finally we have to fit these approximations  to ideas getting rid
of the endomorphisms outside $R$, see Section \ref{constrUT}. We
need the Strong Black Box as discussed and proven in terms of
model theory in Eklof, Mekler \cite[Chapter XIV]{EM}. Note that
this prediction principle is stronger then (Shelah's) General
Black Box, see \cite[Appendix] {CG}. The Strong Black Box is also
restricted to those particular cardinals mentioned in the theorem.
However, here we will apply a version of the Strong Black Box
stated and proven on the grounds of modules in ordinary, naive set
theory, which can be found in  a recent paper by G\"obel, Wallutis
\cite{GW}, see also \cite{GT}. In order to show $\End G = R$
well-known arguments for realizing rings as endomorphism rings
must be modified because the final ring and its action are only
given to us at the very end of the construction:  We will first
replace the layers $G_\a$ from the construction by a new
filtration only depending on the norm. Note that the members of
the new filtration of the right $R$-module $G$ must be right
$R_\a$-submodules for suitable subrings $R_\a$ of $R$ to have
cardinality less than $|G|$. But they are still good enough to
kill unwanted endomorphisms referring to the prediction used
during the construction. Moreover note that the two tasks
indicated in the last two paragraphs must be intertwined and
applied with repetition. While the final $G$ is an $\ale$-free
abelian group, hence torsion-free, it is not hard to see that $G$
is torsion as an $R$-module: If $0\ne g\in G$, then we can choose
distinct elements $g',x,y\in\pp G$ such that $ng' = g,$ and
$x-y\in \pp G$. Hence there are distinct unit elements
$r_x,r_y,r_{xy}\in U(R)=\pm F$ such that $g'r_x =
x,g'r_y=y,g'r_{xy}= x+y$. The endomorphism $r_x+r_y$ does not
belong to $\pm F$, in particular it can not be $r_{xy}$. Hence $r
= r_x+r_y-r_{xy}\ne 0$ but $g'r= x+y-(x+y)=0$ and $g$ is torsion.
It is worth noting that the result can be strengthened under
$V=L$, where we get strongly-$\l$-free groups of cardinality $\l$
as in Theorem \ref{super} for each regular, uncountable cardinal
$\l$ which is not weakly compact. In this case the approximations
in Section \ref{partialauto} can be improved, replacing
`$\ale$-free' by `free' at all obvious places. The main result of
this section will then be a theorem on free groups $G$ with a free
(non-abelian) group $F\subseteq \Aut G$ acting uniquely transitive
on $G$. Also Section \ref{constrUT} must be modified: The Strong
Black Box \ref{sbb} must be replaced by $\Diam$ following
arguments similar to \cite{DG}.

\section{Warming up: Construction of a special group}\label{warmup}

We begin with a particular case of an old theorem and discuss
extra properties of the constructed group. Part of this
proposition will be used in Section \ref{partialauto}.

\begin{proposition}\label{oldgroup} Let $\kappa$ be a cardinal with
$\kappa^{\aln}=\kappa$, $F$ be a free (non-abelian) group of rank
$< \kappa$ and $R =\Z F$ be its integral group ring. Then there is
a group $G$ with the following properties.
\begin{enumerate}
\item $G$ is an $\ale$-free abelian group of rank $\kappa$ with $\End G =
R$.
\item $G$ is torsion-free as an $R$-module.
\item $\Aut G = \pm F$
\item If $\v\in \End G$, then $\v$ is injective.
\end{enumerate}
\end{proposition}

\Proof Note that the integral group ring $R=\Z F$ has free
additive group $R^+$ with basis $F$. We can apply a main theorem
from \cite{CG} showing the existence of an $\ale$-free abelian
group $G$ with $\End G = R$. The free group $F$ is orderable (i.e.
has a linear ordering which is compatible with multiplication by
elements from the right), see Mura, Rhemtulla \cite[p. 37]{MR}.
However note, that torsion-free groups may be non-orderable, see
\cite[pp. 89 - 95, Example 4.3.1.]{MR}. The integral group ring
$\Z F$ of any orderable group $F$ satisfies the unit conjecture,
this is to say that the units of $R=\Z F$ are the obvious ones,
hence $U(\Z F) = \pm F$, see Sehgal \cite[ p. 276, Lemma
45.3]{Se}. Moreover, any group ring which satisfies the unit
conjecture also satisfies the zero divisor conjecture, hence $R$
has no zero-divisors, see Sehgal \cite[p. 276, Lemma 45.2]{Se}.
Therefore $R$ is torsion-free as an $R$-module.

Now the remaining part of the proof is easy: $\Aut G = U(R) = \pm
F$ and if $0\ne g\in G$, then $g\in \oplus R \subseteq G$ because
$G$ is also an $\ale$-free $R$-module by construction in
\cite{CG}. If we consider multiplication of $g$ by some $r\in R$
on a non-trivial component of $g$ in this free direct sum, then
$r=0$ because $R$ is torsion-free as an $R$-module. Hence $G$ is
torsion-free as an $R$-module. Any $\v\in \End G = R$ is scalar
multiplication by a suitable $r\in R$ hence injective because $G$
is a torsion-free $R$-module. \fine

We will use Proposition \ref{oldgroup} in Section
\ref{partialauto}. We get more out of it if we know that a
particular endomorphism is pure:

\begin{lemma}\label{pureimage} Let $F$ be the free (non-abelian) group and
$\End G = \Z F$ be the endomorphism ring of the $\ale$-free
abelian group $G$ given by Proposition \ref{oldgroup}. If $\v \in
\Z F\setminus \pm F$, then $\v$ is a monomorphism and not onto. If
$0 \ne \v \in \Z F$, then $\v$ is pure in $\Z F^+$ if and only if
$\Im \v$ is pure in $G$.
\end{lemma}
\Proof All endomorphisms of $G$ in Proposition \ref{oldgroup} are
monomorphisms as shown there. If $\v\in R = \Z F =\End G$ would be
onto, then $\v$ must be an automorphism, thus $\v \in U(R)$, which
is $\pm F$; and this was excluded.

We come to the last assertion. We shall write $0\ne \v = r \in R$
and suppose that $r = n r' \ (n \ne \pm 1) $ is not a pure element
of $R^+$. Note that $nG \ne G$, hence $Gr'\ne Gnr'$ and we can
pick an element $g\in Gr'\setminus Gr$ which is mapped into $Gr$
under multiplication by $n$. Hence $Gr$ is not pure in $G$.
Conversely let $r$ be pure in $R$ and consider any  $g\in G$ such
that $gp \in Gr$ for some prime $p$. Hence $gp = g'r$ and by
construction of $G$ (just note that $G$ is $\ale$-free as
$R$-module) there is $h\in G$ such that $g'\in hR$ and $hR$ is a
pure subgroup of $G$. Hence also $g\in hR$ and we can write $g =
hr_g, g'= hr_{g'}$ which gives $hpr_g = h r_{g'}r$ and $pr_g =
r_{g'}r$ because $G$ is $R$-torsion-free. Using that $p$ cannot
divide $r$ by purity in $R$ and that $r,r_{g'}$ are elements of
the group ring $\Z F$ we can write $r_{g'} = r' p$ for some $r'\in
R$. Finally $gp = g'r = (h r'p)r$, hence $g = h r'r\in Gr$ and
$Gr$ is pure in~ $G$.~\fine

If we replace \cite{CG} in the proof of Proposition \ref{oldgroup}
by \cite{DG}, then we can strengthen Proposition \ref{oldgroup} in
the constructible universe L. We get a
\begin{corollary}  Let $\kappa$ be a regular,
uncountable cardinal which is not weakly compact such that
$\Diam_\kappa$ holds and let $F$ be a free (non-abelian) group of
rank $< \kappa$ and $R =\Z F$ be its integral group ring. Then
there is a strongly-$\kappa$-free abelian group $G$ of rank
$\kappa$ with $\End G = R$ and properties $(ii),(iii)$ and $(iv)$
of Proposition \ref{oldgroup}.
\end{corollary}
Recall that $G$ is $\kappa$-free if all subgroups of cardinality
$<\kappa$ are free, and $G$ is strongly $\kappa$-free if also any
subgroup of cardinality $<\kappa$ is contained in a subgroup $U$
of cardinality $<\kappa$ such that $G/U$ is $\kappa$-free as well.

\section{Growing partial automorphisms}\label{partialauto}

Besides the set $\pp G$ of pure elements of a group $G$ we
consider a particular subset $\pAut G$ of all partial
automorphisms $\v$ of $G$. Here $\v$ is an isomorphism with domain
$\Dom \v$ and range $\Im \v$ subgroups of $G$. The inverse
isomorphism will be denoted by $\v^{-1}$. However note that
$\v^{-1}$ is not the inverse of $\v$ as a member of $\pAut G$
because $\v \v^{-1} = \v^{-1} \v = 1$ only holds if $\Dom \v = \Im
\v = G$. If we want to stress this point, then we call $\v^{-1}$ a
weak inverse element of $\v$. Surely $0\in \Dom \v$ but it will
happen often that $\v\ne 0$ but $\v^2 = 0$ for partial
automorphisms $\v$. Here we denote with $0$ the trivial partial
automorphism with $\Dom 0 = 0 (=\{0\}) $.

Because we are working exclusively with $\ale$-free groups, we
require that $\v\in \pAut G$ if and only if $\v$ is a partial
automorphism and $G/\Dom \v, G/\Im \v$ are $\ale$-free abelian
groups. The composition of partial automorphisms $\v, \psi $ is
again a partial automorphism with $\Dom(\v \psi) = (\Im \v \cap
\Dom \psi)\v^{-1}$ and range $\Im (\v \psi) = (\Im \v\cap \Dom
\psi)\psi$. Thus products of partial automorphisms of $G$ act
naturally on $G$ as partial automorphisms and domain and range are
well defined. If $\psi, \v \in \pAut G$, then we want to show that
$\psi^{-1}, \v\psi\in\pAut G$. Hence it is enough to check the
freeness condition. If we replace $\psi$ by $\psi^{-1}$, then only
domain and image are interchanged, thus trivially
$\psi^{-1}\in\pAut G$. It remains to consider domain and range of
$\v\psi$. Passing to an inverse, as just noted, it is enough to
deal with $\Dom(\v\psi)$. We already observed that
\begin{eqnarray}\label{product} \Dom( \v\psi) =  (\Im \v \cap
\Dom \psi)\v^{-1}.
\end{eqnarray}
\relax From $\psi \in \pAut G$ follows that $G/\Dom \psi$ is $\ale$-free,
hence
$$\Im \v/(\Im \v \cap \Dom \psi) \cong (\Im \v + \Dom \psi)/\Dom
\psi \subseteq G/\Dom \psi$$ is $\ale$-free. We apply $\v^{-1}$
and (\ref{product}) to see that $\Dom \v/\Dom (\v\psi) $ is
$\ale$-free. Moreover $\v \in \pAut G$, and therefore $G/\Dom \v$
is $\ale$-free, hence  $G/\Dom (\v\psi)$ is $\ale$-free as
desired.

We arrive at our first
\begin{lemma}\label{first}  The set
$\pAut G $ of all partial automorphism $\v$ of $ G$ with $G/\Dom
\v$ and $G/\Im \v$ both $ \ale$-free abelian groups is a submonoid
of all partial automorphisms with $1=\id_G$ and $-1= - \id_G$
acting as multiplication by $1$ and $-1$ respectively, which is
closed under taking (weak) inverses.
\end{lemma}

Moreover, if $\FF \subseteq \pAut G$, then $\FFF \subseteq \pAut
G$ is the submonoid of all products taken from the set $\{ \pm 1
\}\cup \FF \cup \FF^{-1}$, where $\FF^{-1} = \{\psi^{-1} : \psi
\in \FF\}$.

We begin with an observation which allows us to consider induced
partial automorphisms on a factor group.

\begin{observation}\label{inducedpaut} If $U\subseteq G$ are
abelian groups and $\v\in \pAut G$ with $(\Im \v \cap U)\v^{-1}
\subseteq U$ and $(\Dom \v \cap U)\v\subseteq U$, then $\v$
induces a partial automorphism $\vbar$ of $\Gbar$ where $\Gbar
=\{\gbar = g+U: g\in G\}$ taking $\gbar$ to $\overline{g\v}$ for
any $g\in\Dom\v$. Moreover $\Dom\vbar =\overline{\Dom\v}$ and
$\Im\vbar =\overline{\Im\v}$.
\end{observation}

\Proof If $g\in G$ and $\gbar \vbar =\overline{0}$ in $\Gbar$,
then $g\v = g'\in U$ and $$g=g'\v^{-1}\in (\Im \v \cap U)\v^{-1}
\subseteq U,$$ hence $\gbar = \overline{0}$ and $\vbar\in \pAut
\Gbar$. The other assertions are also obvious.\fine

In order to show Proposition \ref{notonto} we relate elements of
free (non-abelian) groups and elements in $\pAut G$. It is
important to be able to work with elements of $\pAut G$ acting on
a partial free basis of $G$. To be precise, we will need the
following definition extending freeness from $G$ to $\pAut G$.

\begin{definition} \label{partfree} Let $\FF = \{\v_t :t\in u\}
\subseteq \pAut G$ be a finite set of partial automorphisms. Then
$(G,\FF)$ is called $\ale$-free if any countable subset of $G$
belongs to a countable subgroup $X\subseteq G$ with basis $B$ and
the following properties for any $\v \in \FFF$.
\begin{enumerate}
\item $G/X$ is $\ale$-free.
\item $\v$ induces a partial injection on $B$, that is, if
$b\in B\cap \Dom\v$, then also $b\v\in B$.
\item $X\cap \Dom \v =\langle B\cap \Dom \v\rangle$.
\end{enumerate}
\end{definition}
Passing to weak inverses, it follows from $(iii)$ that also $X\cap
\Im \v =\langle B\cap \Im \v\rangle$. Moreover $X\cap \Dom \v$ and
$X\cap \Im \v$ are summands of $X$ with free complements generated
by $B\setminus \Dom \v$ and $B\setminus \Im \v$, respectively. It
also follows that $G$ is $\ale$-free. We can ease arguments in
Lemma \ref{morepa} and Lemma \ref{ext} to note here that we only
need a partial basis $b\FFF$ (a subset of $B$) with the property
$(ii)$ and $\langle b\FFF\rangle \cap \Dom \v =\langle b\FFF\cap
\Dom \v\rangle$; see Proposition \ref{notonto}.

Next we relate basis elements of free non-abelian groups and
partial automorphisms with care. Suppose the set $\FF =\{\v_t :
t\in J\}$  generates a free group $\FFF$ and as in Definition
\ref{maindef} there is a map $\pi :\FF \arr \pAut G$ (acting on
the left), then this map can be extended to $\FFF$. The extension
is unique if we restrict ourselves to {\bf reduced elements} $\v
=\v_1\dots \v_n$ in $\FFF$ with $\v_i\in \FF\cup \FF^{-1}$ and
define naturally $\pi(\v) =\pi(\v_1)\dots \pi(\v_n)$. However, if
$\v_1, \v_2\in \FFF$ are reduced and $\v$ is the reduced element
which coincides with the formal product $\v_1\v_2$ in $\FFF$, then
only $\pi(\v_1)\pi(\v_2) \subseteq\pi(\v)$ holds as a graph and
this means $\Dom(\pi(\v_1)\pi(\v_2)) \subseteq\Dom \pi(\v)$ and
$\pi(\v)\restr \Dom(\pi(\v_1)\pi(\v_2)) =\pi(\v_1)\pi(\v_2)$. Thus
we have equality if also the formal product $\v_1\v_2$ is reduced.
\begin{definition}\label{U-property} With $\pi :\FF \arr \pAut G$
as above we say that $\pi$ (or $\pi(\FF)=\{\pi(\v) :\v\in \FF\}$)
satisfies the U-property if $\v =\v'$ for any reduced elements
$\v,\v'\in \FFF$ with $x\pi(\v) =x\pi(\v')$ and some $x\in
\Dom(\v)\cap \Dom(\v')\cap \pp G$.
\end{definition}
The last definition is crucial for this paper because it is the
microscopic version of U-groups discussed in the introduction. We
also must pass from groups $G^\x$ with this U-property to suitable
extension $G^\y$ with the U-property. All this we encode into our
main definition of quintuples $\x$ and their extensions. Normally
our maps will act on the right, but we allow three exceptions, the
maps $\e, \pi$ and $h$ below. Also $\PP_{\aln}(J)$ denotes all
finite subsets of the set $J$.

\begin{definition}\label{maindef} Let $\KK$ be the family of all
quintuples $\x = (G,\FF,\e,\pi,h) =
(G^\x,\FF^\x,\e^\x,\pi^\x,h^\x)$ such that the following holds.
\begin{enumerate}
\item $G$ is an $\ale$-free abelian group.
\item $\FF = \{\v_t : t\in J\}$ is a set of free generators  $\v_t$ indexed
by $J = J^\x$ of a group $\FFF$.
\item $\e :J\to \{1,-1\}$ is a map.
\item $\pi:\FF \arr \pAut G$ is a map which satisfies the
U-property. We shall write $\pi^\x (\v_t) = \v^\x_t$ and omit $\x$
if the meaning is clear from the context.
\item $h:\PP_{\aln}(J) \arr \Im(h)$ is a partial function from
$\Dom h \subseteq \PP_{\aln}(J)$. If $u\in \Dom h$ and $U = h(u)$,
$\FF =\{\v_t : t\in u\}$, then the following conditions must hold.
\begin{itemize}
\item[(a)] $U$ is a countable subgroup of $G$ and
$(\Dom \v^\x \cap U)\v^\x \subseteq U$ for all $\v \in \FFF$;
hence $\v^\x$ induces $\vbar^\x \in \pAut (G/U)$; see Observation
\ref{inducedpaut}. Let $\Gbar = G/U$ and $\FFbar =\{\vbar^\x :
\v\in \FF\}$.
\item[(b)] $(\Gbar,\FFbar)$ is $\ale$-free; see Definition \ref{partfree}.
\end{itemize}
\end{enumerate}
\end{definition}

It follows that $\Gbar$ above is $\ale$-free. From Definition
\ref{maindef} $(iv)$ and Lemma \ref{first} follows
\begin{corollary} If $\x\in \KK$ then $\langle \v_t^\x :t\in
J^\x\rangle\subseteq \pAut G^\x$.\end{corollary}

Hence $G^\x/\Dom\v^\x$ is $\ale$-free for all $\v\in \FFF$. We
will carry on this condition inductively, just checking the
generators in $\FF$ and using the following simple
\begin{testlemma}\label{testlemma} If $U\subseteq G$ are groups and any
countable subset of $G$ is contained in a countable subgroup $X$
with free generators $B_1\cup B_2$ such that $B_1\subseteq U$ and
$U \cap \langle B_2\rangle = 0$ then $G/U$ is
$\ale$-free.\end{testlemma}

The same test lemma will be used inductively for $U$ in Definition
\ref{maindef} $(v)(b)$. We will pass from groups $G^\x$ to larger
groups $G^\y$ related to $\x, \y \in\KK$ by taking pushouts or
unions. This is reflected in the next definition (in particular
condition $(iii)$) of an ordering on $\KK$. This is the final step
before we can start working.

\begin{definition}\label{ordnung} Let $\x\le \y \ (\x,\y\in \KK) $  if
the following holds for $\x = (G^\x,\FF^\x,\e^\x,\pi^\x,h^\x)$ and
$\y = (G^\y,\FF^\y,\e^\y,\pi^\y,h^\y)$.
\begin{enumerate}
\item $G^\x\subseteq G^\y$ and $G^\y/G^\x$ is $\ale$-free.
\item $\pi^\y$ extends $\pi^\x$ in the weak sense $(\pi^\x \preceq
\pi^\y)$, that is $J^\x\subseteq J^\y$ (equivalently
$\FF^\x\subseteq \FF^\y$) and also $\v_t^\x \subseteq \v_t^\y$
extends for all $t\in J^\x$.
\item If $t\in J^\x$, then one of the following cases holds
\begin{itemize}
\item[(a)] $\e^\x(t) = \e^\y (t)$ and $\v^\x_t = \v_t^\y$.
\item[(b)] $G^\x\subseteq \Dom \v_t^\y\cap \Im \v_t^\y$.
\item[(c)] $\e^\x(t)= 1 = -\e^\y(t)$ and $G^\x \subseteq \Dom
\v_t^\y$.
\item[(d)] $\e^\y(t)= 1 = -\e^\x(t)$ and $G^\x \subseteq \Im
\v_t^\y$.
\end{itemize}
\item $h^\x \subseteq h^\y$ extends (i.e., if $h^\x(u)\subseteq G^\x$,
then $h^\x(u)=h^\y(u)\subseteq G^\y$).
\item If $u\in \Dom h^\x$ and $\Gbar^\x =^\df
G^\x/h^\x(u)\subseteq \Gbar^\y =^\df G^\y/h^\x(u)$, then any basis
$B$ of a countable subgroup $X\subseteq \Gbar^\x$ as in Definition
\ref{maindef} extends to a basis $B'$ of some countable subgroup
$X'$ of $\Gbar^\y$ which also satisfies Definition \ref{maindef}.
\end{enumerate}
\end{definition}
\begin{proposition}\label{notonto} Suppose that $\x < \y$ in $\KK$
and $u\in \Dom h^\y, |G^\y|>\aln, G^\y = \Dom \v_t^\y = \Im
\v_t^\y$ for every $t\in u$. If $F$ is the group freely generated
by $\{\v_t : t\in u \}$ and $\theta = \sum_{i\in I} a_i \theta_i
\in \Z F$ is an element of the integral group ring, $a_i \ne 0$
for all $i\in I$ and if the $\theta_i$s are pairwise distinct
(reduced) elements of $F$, such that $\theta^\y =\sum_{i\in I} a_i
\theta_i^\y$ is bijective, then $I$ is a singleton and its
coefficient is $1$ or $-1$.
\end{proposition}
\Proof If $F=\langle \v_t : t\in u \rangle$, then by hypothesis
$F^\y =\langle \v_t^\y : t\in u\rangle$ is a free subgroup of
$\Aut G^\y$. If also $\theta = \sum_{i\in I} a_i \theta_i \in \Z
F$ is as above, then $\theta^\y \in \Z F^\y$ and we may assume
that $\ker \theta^\y =0$. It remains to show that $\theta^\y$ is
surjective if and only if $I =\{0\}$ is a singleton and $a_0= \pm
1$. If $I =\{0\}$, then it is clear that $\theta^\y$ is surjective
if and only if $a_0 =\pm 1$. Hence we may assume that $\theta^\y$
is an isomorphism, and $|I|>1$ for contradiction. In order to
apply Definition \ref{maindef} we pass to the quotient $\Gbar^\y =
G^\y/h(u)$ and to the induced maps $\vbar_t$, which we rename
again as $G^\y, \v_t$. It follows that $h(u)\theta^\y\subseteq
h(u)$ and silently we assume that $h(u)$ is invariant under
$(\theta^\y)^{-1}$; otherwise we must enlarge $h(u)$ by a back and
forth argument such that the quotient satisfies again Definition
\ref{maindef} $(v)$. Also $\Gbar^\y = G^\y/h(u)\ne 0$ because
$|h(u)|=\aln < |G^\y|$.

If $X\ne 0$ is a countable subgroup of $G^\y$, then $X$ is free.
We may assume without restriction that $X\theta^\y \subseteq X,\
XF^\y\subseteq X$ and $X(\theta^\y)^{-1} \subseteq X$ and
$\theta_X=^\df \theta^\y\restr X\in \End X$. If $x\in X$, then $x
= g\theta^\y \in G^\y\theta^\y = G^\y$, thus $g =
x(\theta^\y)^{-1}\in X$ and $\theta_X$ is also surjective (on
$X$). We can start with some $X'$ with a special basis $B'\ne
\emptyset$ as in Definition \ref{maindef} $(v)$ (the weak version
mentioned after the Definition \ref{partfree}) and let $X$ be its
closure as above. Then $X'$ will be a summand of $X$ because
$G/X'$ is $\ale$-free. Hence $B'$ extends to a basis $B$ of $X$:
The maps $\v_X =\v^\y\restr X$, $(\v\in F)$ (by hypothesis) are
total automorphisms of $X$, thus all automorphisms of
$F_X=\{\v_X:\ \v\in F\}$ are permutations of $B'$ when restricted
further to $B'$.

Let $G$ be the group given by Proposition \ref{oldgroup}. Note
that $\v_t\in F \ (t\in u)$ is given by Proposition \ref{notonto}
and $\Z F = \End G$, hence any $\v_t$ can be viewed as an
automorphism of $G$. In order to distinguish it from the element
in $F$ we will call the automorphism $\v_t^*\in \Aut G$ and $F$
becomes $F^*$. The mapping ${}^*$ extends naturally to all of $\Z
F$, (by the identification $\Z F = \End G$), thus $\theta^*=
\sum_{i\in I} a_i \theta_i^* \in \End G$, where $\theta_i^* \in
F^*$. From $|I| > 1$ and Lemma \ref{pureimage} it follows that
$\theta^*$ is not surjective. Let $y\in G\setminus G\theta^*$
which will help us showing that $\theta_X$ can not be surjective
either, in fact we want to show that $B'\cap X\theta_X =
\emptyset$. Fix an element $c\in B'$ and define a map $\Phi:B\to
G$ such that $c\Phi = y$. If $b\in B$ and there is $\v\in F$ such
that $c\v_X = b$, then put $b\Phi = c\Phi \v^*$. If $\v$ exists,
then it is unique by the U-property. Hence $\Phi$ is defined on
$cF_X (=cF^\y)$. If $b\in B'\setminus cF_X$, then let $b\Phi =0$.
Hence $\Phi$ is well-defined on $B$ and extends uniquely to an
homomorphism $\Phi: X\arr G$. It follows $b\v_X \Phi = b\Phi\v^*$,
hence $b\theta_X\Phi = b\Phi \theta^*$ for all $b\in B$, i.e.
commuting with $\Phi$ replaces the ${}_X$ by ${}^*$. If $c\in
X\theta_X$, then there is $x= \sum_{b\in B}x_bb\in X$ with
$x\theta_X = c$. Thus $c = x\theta_X =\sum_{b\in B}x_b(b\theta_X)$
and we apply $\Phi$ to this equation to get the contradiction
$$y =c\Phi =\sum_{b\in B}x_b(b\theta_X)\Phi =\sum_{b\in
B}x_b(b\Phi)\theta^* \in G\theta^*.$$ Thus $\theta_X$, hence
$\theta^\y$ is not surjective.  \fine

It is convenient to check the U-property by the following simple
characterization.

\begin{proposition}\label{U-test} Let $\FFF$ be the group freely
generated by $\FF$ and $\pi: \FF \arr \pAut G^\x$ be a map as in
Definition \ref{maindef} with $\pi(\v_t)= \v_t^\x$ for all $t\in
J$. Then $\FF$ satisfies the U-property if and only if any reduced
product $\v\in \FFF$ with $x\v^\x = x$ for some $x\in
\Dom\v^\x\cap\pp G$ is $\v = 1\in \FFF$. \end{proposition}

\Proof If we can choose a reduced element $\v\in \FFF$ with
$x\v^\x = x = (x1^\x)$ for some $x\in \Dom \v^\x \cap \pp G$, then
$\v = 1$ follows by the U-property of $\FF$. Conversely, if there
are reduced elements $\v,\psi\in \FFF$ with $x\v^\x = x\psi^\x$
for some $x\in \Dom \v^\x\cap \Dom \psi^\x\cap \pp G$, then we can
write $x = x\v^\x (\psi^\x)^{-1}$. Hence $x\in \Dom \v^\x
(\psi^\x)^{-1}$ and we can cancel $\v \psi^{-1}$ to get a reduced
$\theta\in \FFF$ with $\theta =\v \psi^{-1}$ in $\FFF$. From
$x\in\Dom\v^\x (\psi^\x)^{-1}\subseteq \Dom \theta^\x$ it follows
$x = x\theta^\x$. We have $\theta = 1$ by hypothesis, and $\v
=\psi$ follows. \fine

The last proposition shows that the U-property is a strong
restriction on $\pi(\FF)$. If only $x\v^\x = x$ for a reduced $\v$
and pure $x\in G$, then $\v = 1$. However note that if $\v\in
\FFF\setminus \{1\}$, then $x\v^\x (\v^\x)^{-1} = x$ for some
$x\in \pp G$, hence $\v^\x (\v^\x)^{-1}\subseteq \id_G$ but not
$\v^\x (\v^\x)^{-1}= \id_G$ because $\v\v^{-1}$ is not reduced.

The next lemma will be used to make the desired group `more
transitive'. We want to isolate the argument on the existence of
$h(u)$: If $\x = (G,\FF,\e,\pi,h)\in \KK, \ \v_0\in \pAut G$ with
$0\notin J$ as in Lemma \ref{morepa}, then we extend
$h:\PP_{\aln}(J) \arr \Im(h)$ to $h':\PP_{\aln}(J') \arr \Im(h')$
where $J' = J\cup \{0\}$. If $u\in \Dom h$, then $h'(u)=h^\x(u)$.
If $x,y$ and $\x$ are from Lemma \ref{morepa}, then let $h(u')$
for $u' = u\cup \{0\}$ be a countable subgroup $U'\subseteq G$
containing $\{x,y\}\cup h(u)$ such that $\Gbar = G/U'$ and the
induced maps $\FFbar^\y$ satisfy Definition \ref{maindef} $(v)$
(the weak version mentioned after Definition \ref{partfree} will
suffice). Note that we only use extensions of groups $G^\x$ as in
Lemma \ref{morepa} or Lemma \ref{ext} and unions of ascending
chains of such groups. By a back and forth argument, and a moments
reflection about the action of the extended partial automorphisms
by the pushouts, it follows that such a countable subgroup $U'$
exists.

\begin{lemma}\label{morepa} {\bf to get more partial automorphisms.} If
$\x = (G,\FF,\e,\pi,h)\in\KK$ and $x,y\in\pp G$ such that
$x\v^\x\ne \pm y$ for all $\v\in \FFF$, then let $\v_0^\y:x\Z\arr
y\Z (x\arr y)$ be the natural isomorphism. If $\FF^\y = \FF\cup
\{\v_0\}, J^\y = J\cup \{0\}, \e^\y = \e \cup \{(0,1)\}, \pi^\y =
\pi^\x\cup \{(\v_0,\v_0^\y)\}$ and $h^\y =h'$ as above, then $\x<
\y = (G,\FF^\y,\e^\y,\pi^\y,h^\y)\in \KK$
\end{lemma}

\Proof Obviously $\x\le \y$. Also $\v_0^\y \in \pAut G$ because
$x\Z, y\Z$ are pure subgroups of $G$ and $G$ is $\ale$-free, hence
$G/x\Z$ and $G/y\Z$ are $\ale$-free. But it is not clear at the
beginning that $\y$ satisfies the U-property. We will check this
with Proposition \ref{U-test}.

Let be $\FF^\y = \FFs, \v_0 =\eta$ and $\v\in \FFFs$. We write $\v
= \v_1 \eta^{\e_1} \v_2\dots \eta^{\e_{k-1}}\v_k$ with $0 \neq
\e_i\in \Z$ and $\v_i\in \FFF$ reduced and assume that all
$\v_i$'s are different from $\pm 1$, except possibly $\v_1,\v_k$.
Now we assume that $z\v^\y = z$ for some $z \in \Dom \v^\y\cap \pp
G$ and want to show that $\v =1$.

However next we claim, that the product $\v$ must be very special
and show first that $\e_i = \pm 1$ for all $i < k$. If this is not
the case, then some $\eta^2$ or $\eta^{-2}$ is a factor of $\v$.
We may assume that $\eta^2$ appears. Note that $\Dom (\eta^\y)^2 =
(\Im (\eta^\y) \cap \Dom\eta^\y)(\eta^\y)^{-1}$, and $\Im \eta^\y
\cap \Dom\eta^\y = \Z y \cap \Z x = 0$ by the choice of $x,y$.
Hence $(\eta^\y)^2 = 0$ and $\v^\y = 0$ is a contradiction,
because $0\ne z\in \Dom \v^\y$, so the first claim follows. Next
we show that
\begin{eqnarray}\label{short} \v = \pm \v_1\eta^{\e_1} \v_2 .
\end{eqnarray}
We look at the path of $z$, the set $[z]$ of all consecutive
images of $z$: $$z_0 = z , z_1 = z\v_1^\y, z_2 =
z_1(\eta^\y)^{\e_1},\dots, z_{2k-1} = z_{2k-2}\v_k^\y$$

In order to apply $(\eta^\y)^{\e_1}$ to $z_1$ we must have $z_1\in
\Dom (\eta^\y)^{\e_1}$, but $\Dom (\eta^\y)^{\e_1}$ is either $\Z
x$ or $\Z y$, hence $z_1$ is one of the four elements of the set
$V = \{\pm x,\pm y\}$ by purity. Inductively we get $z_i \in V$
for all $0< i < 2k-1$. Suppose that $\eta^{\e_2}$ appears in $\v$,
then $z_3 = z_2\v_2^\y\in V$ because $z_4 = z_3(\eta^\y)^{\e_2}$
is defined and $z_3$ is pure. However $x\in \Dom \v_2^\y$ or $y\in
\Dom \v_2^\y$, respectively. Hence $\v_2$ is multiplication by
$\pm 1$ on $\Z x$ or on  $\Z y$ or $x\v_2 = \pm y$. The last case
was excluded by our hypothesis on $x,y$ and the first two cases
and the U-property would give $\v_2 =\pm 1$ which also was
excluded. Hence $(\ref{short})$ follows.

Our assumption is reduced to $ z = \pm z\v_1^\y(\eta^\y)^{\e_1}
\v_2^\y $ for some $z\in \pp G$. We may replace $\eta$ by
$\eta^{-1}$, hence $\e_1 = 1$ without loss of generality and $z =
\pm z \v_1^\y\eta^\y \v_2^\y$. We consider the path $[z]$ and have
$z_1 = z\v_1^\y = \pm x$ from purity of $z_1\in \Dom \eta^\y$. It
follows $z_2 =\pm y$ and $z_3 = y\v_2^\y =\pm z$ from our
assumption. Thus $y\v_2^\y\v_1^\y = \pm x$ and $\v_2\v_1\in \FFF$,
which contradicts our choice of $x,y$. Hence only $\v = 1$ is
possible and $\y\in\KK$ follows. \fine

The next lemma will increase domain and image of partial
automorphisms, respectively.

\begin{lemma}{\bf growing the partial automorphisms.}
\label{ext} If $\x = (G,\FF,\e,\pi,h)\in\KK$ with $\FF =\{\v_s
:s\in J\}$ and $t\in J$, then there is $\x\le \y=
(G^\y,\FF^\y,\e^\y,\pi^\y,h)\in\KK$ such that the following holds.
\begin{enumerate}
\item $\FF^\x = \FF^\y$, $\e^\x\restr (J\setminus \{t\}) =
\e^\y\restr (J\setminus \{t\})$, $\e^\x(t) = -\e^\y(t)$ and $G^\y
= G_0 + G_1$ is a pushout with $D = G_0\cap G_1$ and $G^\x =
G_0\cong G_1$.
\item
\begin{itemize}
\item[(a)] If $\e^\x(t) = 1$, then $G_0 = \Dom \v_t^\y$,
$G_1 = \Im \v_t ^\y$ and $D =\Dom \v_t^\x$.

\item[(b)] If $\e^\x(t) = -1$, then $G_0 = \Im \v_t^\y$,
$G_1 = \Dom \v_t^\y$ and $D = \Im \v_t^\x$.
\end{itemize}
\end{enumerate}
\end{lemma}

\Proof The set $J$ and $h$ do not change when passing from $\x$ to
$\y$. Thus we consider $\pi^\x$ next and restrict to $\e^\x =1$
(the case $\e^\x =-1$ follows if we replace $\v_t$ by
$\v_t^{-1}$). For the pushout we let $G^\y = (G\times G)/H$ with
$H =\{(x\v, -x) : x\in \Dom\v_t \}$. If we also say that $U_0 =
(U\times 0)+ H/H\subseteq G^\y$ and $U_1= (0\times U)+H/H$ for any
$U\subseteq G$, then in particular $G^\y = G_0 + G_1$ and $D =^\df
G_0 \cap G_1= (\Im \v_t^\x)_0 = (\Dom \v_t^\x)_1$ by the pushout.
Moreover we identify $G_0 =G^\x$, hence $\Dom \v_t^\x = (\Dom
\v_t^\x)_0=^\df D'$ and $D=\Im \v_t^\x$. The canonical map
$$\v_t^\y: G^\y \arr G^\y \ ((x,0) + H \arr (0,x)+ H)$$
extends $\v_t^\x$ because $((x,0) + H)\v_t^\y = (0,x) + H =
(x\v_t^\x,0) + H$ for all $x\in \Dom \v_t^\x$. Clearly $G_0 = \Dom
\v_t^\y$ and $G_1 = \Im \v_t^\y$. Moreover $G^\y/D =G_0/\Dom
\v_t^\x\oplus G_1/\Im \v_t^\x$ is $\ale$-free, hence also $G^\y$
and $G^\y/G^\x\cong G_1/D$ are $\ale$-free, and the maps $\v_s^\x
=\v^\y_s \ (t\ne s\in J)$ remain the same. It follows that
$\pi^\y:\FF^\y\arr \pAut G^\y$. The existence of a partial basis
satisfying Definition \ref{partfree} was discussed before Lemma
\ref{morepa}. So for $\x\le y\in\KK$ we only must check the
U-property for $\FF$ with the new partial automorphisms from
$\v^\y_s\ (s\in J)$ and apply Proposition \ref{U-test}:

Consider a reduced product $\v = \v_1 \v_t^{\delta_1} \v_2\dots
\v_t^{\delta_{n-1}}\v_n, $ where $1\ne\v_i \in \langle \FF
\setminus \{\v_t\}\rangle$ except possibly $\v_1 =\pm 1$ and $\v_n
=1$.

Suppose that $z\v^\y  = z$ for some $z\in \pp G^\y$ and let
$$  z_0 = z, z_1 =z_0\v_1^\y, t_1 =z_1(\v_t^\y)^{\delta_1},
z_2 = t_1\v_2^\y , t_2 =z_2(\v_t^\y)^{\delta_2},\dots, z_n=
t_{n-1}\v_n^\y $$ be the path $[z]$ of $z$. Thus $z_n=z_0$ by
assumption on $\v$. First we note that $\v_i^\y = \v_i^\x$ for all
$i\le n$ with the possible exceptions for $\v_1 = \pm  1$ or $\v_n
= 1$. If also all the $(\v_t^\y)^{\delta_i}$ can be replaced by
$(\v_t^\x)^{\delta_i}$, then $[z]\subseteq G^\x$ and by the
U-property of $\FF$, $z\v^\y = z\v^\x = z$ it follows $\v = 1$. We
will consider the two cases $z_0 \in G_0$ and $z_0 \in
G_1\setminus G_0$.

First we reduce the second case $z_0 \in G_1\setminus G_0$ to the
first case. Since $z_0\in\Dom \v_1$ it follows $\v_1 = \pm 1$,
hence $z_1 = \pm z_0\in G_1\setminus G_0$. From $z_1 \in \Dom
(\v_t^\y)^{\delta_1}$ it follows $\delta_1 \le -1$ and $t'_1 =^\df
z_1(\v_t^\y)^{-1}\in G_0$. From $z_n = z_0$ it follows
$z_{n-1}(\v_t^\y)^{\delta_{n-1}}\v^\y_n = z_0\in G_1\setminus G_0$
and therefore $\v_n =1$ and $\delta_{n-1} \ge 1$. The equation
$z_0\v^\y = z_0$ reduces to
$$\pm t'_1(\v_t^\y)^{\delta_1+1}\v_2^\x (\v_t^\y)^{\delta_2} \cdots
\v_{n-1}^\x = t'_1 \mbox{ with } t'_1\in G_0,$$ which is the first
case for a new $z =t'_1\in \pp G_0$.

If $z_0 \in G_0$, then also $z_0\in \Dom \v_1^\x$ and $z_1 = z_0
\v_1^\x\in G_0$. We will continue along the path step by step
having two subcases each time, but one of them will lead to a
contradiction. In the first step either $z_1\in D'$ or $z_1\notin
D'$. In any case $\delta_1 =1$ and in the second case $t_1
=z_1\v_t^\y\in G_1\setminus D$, but $t_1\in \Dom \v_2^\y$ so
necessarily $\v_2 =1$ and $n=2$. We get $t_1\v_2^\y = z_2 =t_1$,
and $t_1 = z_0\in G_0$ contradicting $t_1\in G_1\setminus D$. We
arrive at the other case $z_1 \in D'$. Hence $t_1 = z_1\v_t^\x \in
D$ and we must have $t_1\in \Dom \v_2^\x$. Therefore also $z_2
=t_1 \v_2^\x\in G_0$, and $\delta_2 = 1$. Again we have two cases
$z_2\in D'$ and $z_2\notin D'$ with $\delta_2 = 1$, where the
second case leads to a contradiction. Hence $t_2=z_2\v_t^\x\in D$
and we continue until we reach $n$ and all the $\v_t^\y$s are
replaced by the $\v_t^\x$s. Now our first remark applies and $\v =
1$ follows from the U-property for $\FF$. \fine

\begin{lemma}\label{chain} Let $\a$ be a limit ordinal. Then any increasing
continuous chain $\x^j = (G^j,\FF^j,\e^j,\pi^j,h^j)_{ j\in\a}$ in
$(\KK,<)$ obtained by applications of Lemma \ref{morepa} and Lemma
\ref{ext} has a natural supremum $\x = (G,\FF,\e,\pi,h)$ in $\KK$,
where $G =\bigcup_{j\in \a} G^j$, $\FF =\bigcup_{j\in \a} \FF^j$,
$J = \bigcup_{j\in \a} J^j$ and $\pi,h$ are defined below.
\end{lemma}

\Proof The map $h$ extends uniquely, because $h$ is defined on
finite subsets of $J$. Similarly we can handle $\pi$. We define
$\pi: \FF \arr \pAut G$ by taking unions: If $t\in J$, then $t\in
J^j$ for all $i < j\in \a$ and $i\in\a$ large enough. Therefore we
can let $\v_t^\x =^\df \bigcup_{i < j\in \a} \pi^j(\v_t^j)$, and
$\Dom \v_t^\x = \bigcup_{i < j\in \a} \Dom \pi^j(\v_t^j)$. We
found a well-defined partial automorphism $\v_t^\x :\Dom \v_t^\x
\arr G$ and also want that $\pi(\v_t) = \v_t^\x \in \pAut G$.
Hence we must show that $G/\Dom \v_t^\x $ and $G/\Im \v_t^\x $ are
$\ale$-free. Passing to an inverse it is enough to consider
$G/\Dom \v_t^\x $. Either there is a strictly increasing chain
$j_i \in \a \ (i\in I)$ cofinal to $\a$ such that
$\pi^{j_i}(\v_t^{j_i}) \ne \pi^{j_i+1}(\v_t^{j_i+1})$ for all
$i\in I $ or the sequence $\pi^{j}(\v_t^{j})$ becomes stationary,
say at $j_0\in \a$. In the first case Lemma \ref{ext} applies and
$\Dom\pi^{j_i+1}(\v_t^{j_i+1})= G_0^{\x_{j_i}}$ and $\bigcup_{i\in
I} G_0^{\x_{j_i}}= G^\x$ by cofinality, hence $\Dom \v_t^\x = G$
and $G/\Dom \v_t^\x = 0$ is trivially $\ale$-free. In the other
case we have $\Dom\v^\x_t = \Dom\v^{j_0}_t$, hence
$G^{j_0}/\Dom\v^\x_t$ is $\ale$-free by induction. Moreover
$G/G^{j_0}$ is $\ale$-free because $\ale$-free is of finite
character by Pontryagin's theorem (speaking about subgroups of
finite rank). Hence $G/\Dom\v^\x_t$ is $\ale$-free, as required.
Now it is easy to see that $\x\in \KK$: Definition \ref{partfree}
is easily verified by taking union of partial bases. The
U-property carries over to limit ordinals. \fine

The next main result of this section follows by iterated
application of the Lemma~\ref{morepa}, Lemma~\ref{ext} and
Lemma~\ref{chain}. Without danger we now can identify the free
groups $F'$ and $\FFFs$ by the isomorphism $\pi'$ in the theorem.

\begin{theorem} \label{Auto} If $\x = (G,\FF,\e,\pi,h)$ is in
$\KK$, then we can also find $\x' = (G',\FF',\e',\pi',h')$ in
$\KK$ with $\x\le \x'$ such that the following holds.

\begin{enumerate}
\item[(a)] $\FF'$ is a set of automorphisms of $G'$ which freely
generates $\FFFs\subseteq \Aut G$, hence $G'$ is a module over
$R'= \Z \FFFs$.
\item[(b)] $\FFFs$ acts uniquely transitive on $\pp G'$.
\item[(c)]  $|G'|+|\FF'|= |G|+ |\FF|$
\end{enumerate}
\end{theorem}

\Proof We will proceed by induction to get a chain $\x\le \x_n \in
\KK \ (n\in \o) $ taking several steps each time. Finally $\x'$
will be the supremum of the $\x_n$s.

In the first step we apply $|G_0|$-times Lemma \ref{morepa} and
Lemma \ref{chain} to $\x =\x_0$ taking care of all appropriate
pairs of elements in $\pp G_0$ and let $\x_0\le \x_1 =
(G_1,\FF_1,\e_1,\pi_1,h_1)$ be the union of this chain such that
$\langle\FF_1\rangle$ acts transitive on $\pp G_0$. In this case
$G_0 =G_1$ but the other parameters in $\x_0$ increase (along a
chain) by Lemma \ref{morepa}. In the next step we apply
$|G_1|$-times Lemma \ref{ext} and Lemma \ref{chain}  to get
$\x_1\le \x_2= (G_2,\FF_2,\e_2,\pi_2,h_2)$ such that $G_0\subseteq
\Dom \v_t^{\x_2}$ and $G_0\subseteq \Im \v_t^{\x_2}$ for all $t\in
J^{\x_2}= J^{\x_1}$. We continue this way for each $n$. From Lemma
\ref{chain} follows $\x \le \x'\in\KK$, and by construction
$\FFFs\subseteq \Aut G'$ acts uniquely transitive on $\pp G'$ such
that also $(c)$ of Theorem \ref{Auto} holds.  \fine

Now we can restrict elements in $\KK$ to triples and say that $\x
= (G^\x,\FF^\x,\pi^\x)$ belongs to $\KK^*$ if and only if $G^\x$
is an $\ale$-free abelian group, $\FF =\{\v_t : t\in J^\x\}$
freely generates a (non-abelian) group $\FFF$ and $\pi^\x: \FF
\arr \Aut G^\x$ is an injective map such that $ \pi^\x(\FF) = \{
\pi^\x(\v_t) = \v_t^\x : t\in J^\x\}$ satisfies the U-property.
Still we can view $\KK^*$ as a subset of $\KK$ and have an induced
ordering, compare Definition \ref{ordnung}: We have $\x\le \y \
(\x,\y\in \KK^*)$ if $(v)$ and the following holds.
\begin{enumerate}
\item $G^\x\subseteq G^\y$ and $G^\y/G^\x$ is $\ale$-free.
\item $\pi^\x \subseteq \pi^\y$.
\end{enumerate}
The other conditions in Definition \ref{ordnung} are vacuous.

\section{Construction of uniquely transitive groups}\label{constrUT}
In this section we will sharpen Theorem \ref{Auto} which shows
that a group ring $R = \Z F$ for some free group $F$ can be
represented as $R\subseteq \End G$ of some $\ale$-free group $G$
(hence $G$ is an $R$-module), such that the units $U(R) = \pm F$
act uniquely transitive on $G$. If $\pm F\subseteq \Aut G $ is not
necessarily transitive but satisfies the uniqueness property ((
$\v, \v'\in \pm F$ and $\exists x\in \pp G, x\v = x\v')
\Rightarrow \v=\v'$), then we will also say that the pair $(G,F)$
has the U-property. We will modify $G$ such that equality holds,
that is $R = \End G$. Thus we have to kill unwanted endomorphisms,
which is done using the Strong Black Box \ref{sbb}. We need some
preparation to work with this prediction principle.

First we need an $\bS$-adic topology. Let $\bS = \{ q_n : n\in\o
\}$ be an enumeration of a multiplicatively closed set generated
by $1$ and at least one more natural number different from $1$ and
define $s_0 = 1, s_{n+1} = s_n\cdot q_n$ for all $n\in\o$ which
obviously defines a Hausdorff $\bS$-topology on the groups $G$
under consideration. Let $\widehat G$ be the $\bS$-adic completion
of $G$ and $G\pure \widehat G$ naturally, where purity ``$\pure$''
is $\S$-purity.

Next we must formulate the Strong Black Box and adjust our
notations; see Strong Black Box \ref{sbb}. We rely on the version
adjusted to and proven for modules using only naive set theory as
explained in the introduction, see \cite{GW,GT}. Thus we have to
fix a few parameters next. To do so we choose an enumeration by
ordinals $\a\in \l$, with $\l = \mu^+$ a successor cardinal such
that $\mu^{\aln} = \mu$: Let $F_\a$ be a free non-abelian group
with a set of $\mu$ free generators $\FF_\a$, the $\FF_\a$s
constitute a strictly increasing, continuous chain in $\a <\l$. We
will write $R_\a = \Z F_\a$ with $R =\bigcup_{\a\in\l} R_\a$ and
note that $R_\a^+$ and $R^+$ are free abelian groups. By Theorem
\ref{Auto} there is an $R_\a$-module $G_\a$ of cardinality $\mu$
which is cyclic as $R_\a$-module and $\ale$-free as abelian group
and $F_\a$ acts sharply transitive on $\pp G_\a$. We can choose a
free abelian and $\bS$-dense subgroup $B_\a\pure G_\a\pure
\cp{B_\a}$. Also the $B_\a$s constitute a strictly increasing,
continuous chain in $\a <\l$. Passing to $\a +1$ we can let
$B_{\a+1} = B_\a \oplus A_\a$ with $A_\a = \dsum_{i<\rho}\Z a_i$
(and $\rho =\mu$). It helps (when using \cite{GW,GT}) not to
identify $\rho$ and $\mu$ in the formulation of the Strong Black
Box \ref{sbb}.

Thus the construction is based on a free abelian group $B$ of rank
$\l$ and its $\bS$-adic completion $\cp{B}$, in fact the desired
group $G$ will be sandwiched as $B \pure G\pure \cp{B}$. Put $B =
\dsum_{\a<\l} e_\a A_\a$. Then, writing $e_{\a,i}$ for $e_\a a_i
$, we have $B=\dsum_{(\a,i)\,\in\,\l\times\rho}\Z e_{\a,i}$. For
later use we put the lexicographic ordering on $\l\times\rho$;
since $\rho,\l$ are ordinals $\l\times\rho$ is well ordered.
\\ 
For any $g=(g_{\a,i} e_{\a,i})_{(\a,i)\, \in\,\l\times\rho}
\in\cp{B}\subseteq \prod_{(\a,i)\,\in\,\l\times\rho}
\cp{\Z}e_{\a,i}$ we define the {\em support} of $g$ by
$\supp{g}=\{(\a,i)\in\l\times\rho\, |\, g_{\a,i}\ne 0\} $ and the
support of $H\subseteq \cp{B}$ by $\supp{H}=\union_{g\in
H}\supp{g}$; note $\card{\supp{g}}\le \aleph_0$ for all $g\in
\cp{B}$. Moreover, we define the {\em $\l$-support} of $g$ by
$\suppl{g}=\{\a\in\l\,|\, \exists\,i\in\rho:\, (\a,i)\in\supp{g}\}
\subseteq\l$ and the {\em $A$-support} of $g$ by $\suppa{g}=
\{i\in\rho\,|\, \exists\,\a\in\l:\, (\a,i)\in\supp{g}\}
\subseteq\rho\subseteq\l$. Recall that  $e_{\a,i}=a_i e_\a$, where
$a_i\in A_\a$, which explains the use of the notion \lq\lq
$A$-support\rq\rq.

Next we define a {\em norm} on $\l$, respectively on $\cp{B}$, by
$\norm{\{\a\}}=\a+1$ \mbox{$(\a\in\l),$} $\norm{H}=\sup_{\a\in
H}\norm{\a}\ (H\subseteq\l)$, hence $\norm{\a}=\a$, and
$\norm{g}=\norm{\suppl{g}}\ (g\in\cp{B})$, i.e. $\norm{g}=
\min\{\beta\in\l\, |\, \suppl{g}\subseteq\beta\}$. Note,
$\suppl{g}\subseteq\beta$ holds if and only if $g\in\cp{B_\beta}$
for $B_\beta= \dsum_{\a<\beta} e_\a  A_\a$. We also define an {\em
$A$-norm} of $g$ by $\norma{g}=\norm{\suppa{g}}$. Also let
$\l^o:=\{\a<\l\,|\, \cf(\a)=\o \}$.

Finally, we need to define canonical homomorphisms used in the
Strong Black Box for predictions. For this we fix bijections
$g_\gamma\!:\mu\arr\g$ for all $\gamma$ with $\mu\le\g<\l$ where
we put $g_\mu=\id_\mu$ and so $\card{\g}=\card{\mu}=\mu$ for all
such $\g$'s. For technical reasons we also put $g_\g=g_\mu$ for
$\g<\mu$.

\begin{definition}\label{defcan1}
Let the bijections $g_\g\ (\g<\l)$ be as above and put
$\g_{\a,i}=\g_\a\times \g_i$ for all $(\a,i)\,\in\, \l\times\rho$.
We define $P$ to be a {\em canonical summand} of $B$ if $P=\dsum_{
(\a,i)\,\in\, I}\Z e_{\a,i}$ for some $I\subseteq\l\times\rho$
with $\card{I}\le\aln$ such that:
\begin{itemize}
\item
if $(\a,i)\,\in I$, then $(i,i)\in I$; if $(\a,i)\,\in I,
\a\in\rho$ then $(i,\a)\in I$ and
\item
if $(\a,i)\in I$, then $\left(I\cap(\mu\times\mu)\right)g_{\a,i}=
I\cap \a\times i$. \end{itemize} Accordingly, a homomorphism
$\phi\!:P\rightarrow \cp{B}$ is a {\em canonical homomorphism} if
$P$ is a canonical summand of $B$ and $\Im\phi\subseteq \cp{P}$;
we put $\supp{\phi}=\supp{P}$, $\suppl{\phi}=\suppl{P}$ and
$\norm{\phi}=\norm{P}$.
\end{definition}

Note, by the above definition, a canonical summand $P$ satisfies
$\norma{P}\le\norm{P}$. Let $\C$ denote the set of all canonical
homomorphisms. From assumption $\mu^{\aln}=\mu$ follows
$\card{\C}=\l$. We are now ready to formulate the Strong Black
Box:

\begin{bbtheo}\label{sbb}\quad
Let $\l = \mu^+$ and $\mu =\mu^{\aln}$ be as before and let
$E\subseteq\l^o$ be a stationary subset of $\l$.\\ 
Then there exists a family $\C^\ast$ of canonical homomorphisms
with the following properties:
\begin{itemize}
\item[{\bf(1)}]
If $\phi\in\C^\ast$, then $\norm{\phi}\in E$.
\item[{\bf(2)}]
If $\phi,\,\phi'$ are two different elements of $\C^\ast$ of the
same norm $\a$, then\\ 
$\norm{\suppl{\phi}\,\cap\,\suppl{\phi'}}<\a$.
\item[{\bf(3)}]
{\sc Prediction:} For any homomorphism
$\,\psi\!:B\rightarrow\cp{B}\,$ and for any subset $I$ of
$\l\times\rho$ with $\card{I}\le\aln$ the set\\ 
\parbox[t]{12cm}{\centerline{$\{\a\in E\,|\, \exists\, \phi\in\C^\ast:
\norm{\phi}=\a,\,
\phi\subseteq\psi,\,I\subseteq\supp{\phi}\}$}}\\
is stationary.
\end{itemize}
\end{bbtheo}
We will enumerate $\C^*$ taking care of the order of the norm and
can write $\C^* = \{\phi_a: \a < \l\}$ such that
$\norm{\suppl{\phi_\a}} \le \norm{ \suppl{\phi_\b}} $ for all  $\a
<\b<\l$. Note that we distinguish between $\v\in \FF$ and $\phi$
given by the Strong Black Box \ref{sbb}.

The enumeration is now used to find a continuous, ascending chain
$G_\a \ (\a < \l)$ such that $B_\a\pure G_\a \pure \cp{B_\a}$ and
the corresponding $\x_\a = (G_\a, \FF_\a, \pi_\a) \in \KK^*$. Let
$E \subseteq \l^o$ be a fixed stationary subset of $\l^o$ and
choose $\x_0 = (G_0, \FF_0,\pi_0)$ to be any triple given by
Theorem \ref{Auto} for cardinality $\mu$. For any $\b <\l$, let
$P_\b =\Dom\phi_\b$ and suppose that the triples $\x_\b = (G_\b,
\FF_\b,\pi_\b)$ are constructed for all $\b < \a$. If $\a$ is a
limit ordinal, then by continuity we let $\x_\a = \sup_{\b <\a}
\x_\b$, thus $\x_\a = (G_\a, \FF_\a,\pi_\a)$ is well defined and
belongs to $\KK^*$ by Proposition \ref{Uprop} and Theorem
\ref{Auto}; in particular $G_\a = \bigcup_{\b <\a} G_\b$ is an
$R_\a$-module over the integral group ring $R_\a = \bigcup_{\b
<\a} R_\a = \Z F_\a$ of the free group $F_\a = \bigcup_{\b <\a}
F_\b$ and $B_\a\pure G_\a \pure \cp{B_\a}$.

We may assume that $\a = \b +1$ and distinguish two cases, the
trivial situation and the ordinal $\a$ where we have to work.

The first case arrives if we do not want to work or cannot work:
If $\b$ is also a discrete ordinal, then we apply two steps. First
let $G'_\a = G_\b \oplus G^1_\b $ where $G^1_\b= e_\b R_\b$ is a
``new" free summand. In the second step we apply Theorem
\ref{Auto} to get $G'_\a\pure G_\a$ and $\x_\a =(G_\a,
\FF_\a,\pi_\a) \in \KK^*$. Note that in particular $G_\a$ is an
$R_\a$-module with $R_\b\subseteq R_\a$ as integral group rings
(because $F_\b\subseteq F_\a$) and $G_\a/G_\b$ is $\ale$-free.
Using $A_\b$, we can arrange that $B_{\b+1}\subseteq_* G_{\b+1}
\subseteq_* \cp{B_{\b+1}}$. If $\b$ is a limit ordinal and not in
$E$, we proceed similarly. If $\b \in E$, then we also apply the
trivial extension if $\phi_\b$ is scalar multiplication by some
$r\in R_\b$ when restricted to $P_\b$ or if $\Im
\phi_\b\not\subseteq G_\b$. Otherwise we will meet the condition
of the Step Lemma \ref{step} and must work:

Suppose that $\a=\b + 1$ with $\b\in E$  and $\Im \phi_\b
\subseteq G_\b$, $\phi_{\beta} \not\in R_{\b}$. In this case we
try to `kill' the undesired homomorphism $\phi_{\beta}$ which
comes from the black box prediction. (However note that it could
be that $\phi_\b\in R$ later on, so in this case $\phi_\b$ is a
good candidate which should survive the massacre. This we must and
will see clearly at the marked place near the end of the proof!)
Recall that $\norm{\phi_{\b}} \in \l^o$, hence there are $(\b_n,
i_n) \in [\phi_\b ]$ $(n \in \o$) such that $\b_0 < \b_1 < \cdots
< \b_n < \cdots$ and $\sup_{n \in \o}\b_n = \norm{\phi_\b}$.
Without loss of generality we may assume that $\b_n \notin E$ for
all $n \in \o$ and hence $G_{\b_n +1}' = G_{\b_n} \oplus
e_{\b_n}R_{\b_n}$. We put $I= \{(\b_n, i_n) | n < \o \}$. Then
$I_{\l} \cap [g]_{\l}$ is finite for all $g \in G_{\beta}$. We
apply the Step Lemma~\ref{step} to $I, P=\Dom \phi_{\b}$ and
$H=G_\b$. Therefore there exists an extension $\x_\a = \x_{\b+1}$
of $\x_\b$ and an element $y_\b \in G_\a$ such that $y_\b \phi_\b
\notin G_\a$ and $\norm{y_\b} = \norm{\phi_\b} = \norm{P_\b}$.
Thus the chain $G_\a \ (\a<\l)$ is constructed up to the used Step
Lemma \ref{step}. Finally let $G=\bigcup_{\a\in \l}G_\a$.

It remains to show $R =\End G$, the Step Lemma \ref{step} and the
following proposition which ensures that the U-property can be
extended when applying the step lemma. We will use obvious
simplified notations:

\begin{proposition}\label{Uprop} Let $G = \bigcup_{n\in \o} G_n$,
$x_n\in G^1_n$ be pure ($R_n$-torsion-free) elements, where
$G'_{n+1} = G_n \oplus G^1_n\pure G_{n+1}$ and $G_{n+1}$ is
obtained from $G'_{n+1}$ by Theorem \ref{Auto}. If $b\in \cp{G_0},
x = \sum_{n\in\o} x_n s_{n-1} + b$ and $G' = \langle G ,x R_\o
\rangle_*$, then $U(R_\o) = \pm F $ for the free group
$F=\bigcup_{n\in \o} F_n $, the pair $(G',F)$ satisfies the
U-property and $G'$ is $\ale$-free.\end{proposition}

\Proof It is well-known that $G'$ is $\ale$-free, see \cite{CG} or
\cite{EM,GT} for instance. We want to show that the pair $(G',F)$
satisfies the U-property.

If $r\in R_\o$ and $\v\in F$, then $r\in R_n$  and $\v\in F_n$ for
some $n\in \o$ large enough. Hence there is no restriction to
write any $0 \ne y\in G'$ as $y =x^kr + g$ with $x^k = \sum_{n
> k} \frac{ s_{n-1}}{ s_k } x_n + b^k$, $r\in R_n$, $b-b^k\in
G_0$, $g = g_0 + g_1$, $g_0\in \cp{G_0},\ 0\ne g_1\in G_n$ and
$[G_0]\cap [g_1] = [x^k r]\cap [g] = \emptyset$.

Suppose $y\v = y$, hence $(x^kr)\v + g\v = x^kr + g$ and by
continuity $$ \sum_{m > k} \frac{ s_{m-1}}{ s_k } (x_mr)\v + g\v =
\sum_{m > k} \frac{ s_{m-1}}{ s_k } (x_m r) + g.$$ We consider any
$m> n,k$ and note that $x_m$ is a torsion-free element of the
$R_n$-module $G_m$ and also $\v\in U(R_n)$, (hence $G_m$ is
invariant under application of $\v, r $) and $0\ne x_mr\v,
x_m\v\in G_m$. Restricting the support $[x^kr]$ of the last
displayed equation to $x_m$ shows that $\frac{ s_{m-1}}{ s_k }
(x_m r)= \frac{ s_{m-1}}{ s_k } (x_mr)\v$. Hence $0\ne (x_n
r)=(x_n r)\v$ and passing to $h\in \pp G_m$ with $hq = x_n r$ we
get $h\v =h$, hence $\v = 1$ by the U-property for $(G_m,F_m)$. So
$\v$ acts as $\id_{G_m}$ for all $m\in\o$ large enough. Thus $\v$
is the identity on $G$, and the same holds for $G'$ by continuity
($G$ is dense in $G'$ in the $\bS$-topology). \fine

Next we prove the step lemma.

\begin{steplemma}
\label{step} Let $P=\bigoplus\limits_{(\a,i) \in I^*} \Z e_{\a,i}$
for some $I^* \subseteq \l \times \rho$ and let $H$ be a subgroup
of $\B$ with $P \subseteq_* H \subseteq_* \B$ which is $\ale$-free
and an $R$-module, where $R=^\df R^\x = \bigcup\limits_{\a \in
I'}R_{\a}$ with $I' =^\df \suppl{I} = \{ \a < \l : \exists i <
\rho, (\a,i) \in I^* \}$. Assume that $\x =(H,\FF,\pi) \in \KK^*$
(thus $R^\x =\Z \FFF$). Also suppose that there is a set $I = \{
(\a_n, i_n) : n < \o \} \subseteq [P]=I^*$ such that $\a_0 < \a_1
< \cdots < \a_n < \cdots (n < \o)$ are discrete ordinals and
\begin{enumerate}
\item $H=\bigcup\limits_{n < \o}G_{\a_n}$,
$R=\bigcup\limits_{n < \o}R_{\a_n}$; where $R_{\a_n} = \Z \langle
\FF_{\a_n}\rangle$, $\x_n = (G_{\a_n}, \FF_{\a_n}, \pi_n)\in
\KK^*$
\item $I_{\l} \cap [g]_{\l}$ is finite for all
$g \in H$ $(I_{\l}=[I]_{\l})$.
\end{enumerate}
If $\phi : P \rightarrow H$ is a homomorphism which is not
multiplication by an element from $R$, then there exists an
element $y \in \widehat{P}$ and $\x\le \x'=(H', \FF', \pi') \in
\KK^*$ such that $H \subseteq_* H' \subseteq_* \B$, $y \in H'$ and
$y\phi \not\in H'$.
\end{steplemma}

\Proof By assumption $H$ is an $R$-module and hence the
$\bS$-completion $\widehat{H}$ is an $\widehat{R}$-module. Let
$\widehat{\Z}$ be the $\bS$-completion of $\Z$ and let
$x=\sum\limits_{n \in \o}s_n e_{ \a_n, i_n }\in \widehat{P}$. We
will use $y =x\in \widehat{P}$ or $ y = x + b\in \widehat{P}$ for
some $b\in \widehat{P}$ to construct two new groups $H \subseteq_*
H_y$:

Let $H'_y=\langle H, y R \rangle _* \subseteq \B$. By Proposition
\ref{Uprop} and Theorem \ref{Auto} we obtain $\x_y =
(H_y,\FF_y,\pi_y)\in \KK^*$ such that $\x\le \x_y$ and $y\in H_y$.
Clearly $y\phi \in \widehat{H}$. Put $z=y\phi$ and $y=x$. If $z
\not\in H_x$, then we choose $H'=H_x$ and have $\x_x \in \KK^*$
with $\phi \not\in \End H_x$. Also $R_x =\Z\langle \FF_x\rangle$.

If $z \in H_x$, then also $z\in H'_x$ by an easy limit argument
and there exist integers $k$ and $n$ such that $s_kx\phi=g r_n + x
r_n'$ for some $r_n, r_n' \in R$ and $g \in H$. It follows
$x(s_k\phi -r_n') = gr_n$ and $r_n,r'_n\in R_{n'}$ for some
$n'\in\o$. Since $\left(s_k \phi - r_n'\right) \not=0$ there is
$b' \in P$ such that $b=^\df b'(s_k\phi -r_n') \not=0$. Note that
$b'$ has finite support. Moreover, by the cotorsion-freeness of
$H$ there exists $\pi \in \widehat{\Z}$ such that $\pi b \not\in
H$. Let $y= x+ \pi b'$. We claim that $\phi \not\in \End (H_y)$.
By way of contradiction assume that $ s_l(x+\pi b')\phi = g^*
r^*_m + (x+ \pi b')r_m^{**}$  for some integer $l \geq k$ and
elements $r_m^*, r_m^{**} \in R_{m'}$ for some $m'\in\o$ and $g^*
\in H$. Without loss of generality, we may assume $n=m$, hence $
s_l(x+\pi b')\phi = g^* r^*_n + (x+ \pi b')r_n^{**}$ and if
$s=s_l/s_k$, then
\[ s_l(\pi b')\phi = s_l(x + \pi b')\phi -
ss_kx\phi = g^* r_n^* + (x+ \pi b')r_n^{**} - s(gr_n + xr_n') \]
\[= (g^*r_n^* - sgr_n) + \pi b' r_n^{**} +
x(r_n^{**}-sr_n'). \] Since $[\pi b']= [b'], [\pi b'\phi ] = [
b'\phi]$ and $g^*, g \in H$ an easy support argument shows that
$r_n^{**}=sr_n'$, hence $ ss_k (\pi b')\phi = (g^* r_n^* - sgr_n)
+ s\pi b' r_n'  $ and thus $ s\pi  (s_k b'\phi -  b' r_n')= (g^*
r_n^* - sg r_n ) \in H$. By purity we get $\pi(s_k b'\phi - b'
r_n' ) =\pi b \in H$ - a contradiction.  \fine

{\sl Proof of Theorem \ref{super}.\quad} If $R = \End G$, then
also $\Aut G = \pm F$, because $U(R) = \pm F$ and any $r\in
R\setminus U(R)$ viewed as endomorphism of $G$ is not bijective by
Proposition \ref{notonto}. Hence it remains to show that $\End
G\subseteq R$.

First we choose a new filtration of $G$ and define
$$G^\a=^\df \{g\in G, \,\norm{g}<\a \} \  (\a<\l).$$

\begin{lemma}\label{desc2} The new filtration on $G$ has the following
properties.
\begin{itemize}
\item[(a)] $G \cap \widehat{P_\b}\subseteq G_{\b+1}$ for all $\b<\l$;
\item[(b)] $\{G^\a,\  \a<\l\}$ is a $\l$-filtration of $G$;
\item[(c)] If $\a, \b<\l $ are ordinals such that $\norm{\phi_\b}=\a$
then $G^\a\subseteq G_\b$.
\end{itemize}
\end{lemma}

\Proof Inspection of the definitions, or \cite{GW, GT}. \fine

Assume that $\phi \in \End G \backslash R$ and let $\phi'=\phi
\restr B$, hence $\phi' \not\in R$. Let $I=\{(\a_n,i_n)|\ n < \o
\} \subseteq \l\times \rho$ such that $\a_0 < \a_1 < \cdots \a_n <
\cdots$ is a sequence of discrete ordinals and $I_\l \cap [g]_\l$
is finite for all $g \in G$. Note that the existence of $I$ can be
easily arranged, e.g. let $E \subsetneq \l^o, \a \in \l^o
\backslash E, i_n \in \rho \ (n < \o)$ arbitrary and $(\a_n)_{n <
\o}$ any discrete ladder on $\a$.

By Lemma \ref{step} there exist an element $y \in \B$,
$(G,\FF,\pi)\le (G',\FF',\pi')\in \KK^*$ such that $y\in G'$ and $
y \phi' \not\in G'$. (The last innocent sentence uses $\phi'\notin
R$ and not just $\phi'\notin R_\b$ for a suitable $\b$, see our
comment in the construction of $G$.)  By the Black Box Theorem
\ref{sbb} the set
\[ E'=\{ \a \in E | \exists \b < \l : \norm{\phi_\b} = \a,
\phi_\b \subseteq \phi', [y] \subseteq [\phi_\b] \} \] is
stationary since $|[y]| \le \aln$. Note, $[y]\subseteq [\phi_\b]$
implies that $y \in \widehat{P_\b}$. Moreover, the set $C=\{ \a <
\l : G^\a\phi \subseteq G^\a \}$ is a cub in $\l$, hence $E' \cap
C \not= \emptyset$. Let $\a \in E' \cap C$. Then $G^\a \phi'
\subseteq G^\a$ and there exists an ordinal $\b < \l$ such that
$\norm{\phi_\b} =\a$, $\phi_\b \subseteq \phi$ and $y \in
\widehat{P_\b}$. The first property implies that $G^\a \subseteq
G_\b$ by Lemma \ref{desc2} and the latter properties imply that
$\phi_\b \not\in R$. Moreover, $P_\b \subseteq B$ with
$\norm{P_\b}=\a$ and hence $P_\b$, and also $P_\b\phi$ are
contained in $G^\a \subseteq G_\b$.

Therefore $\phi_\b : P_\b \arr G_\b$ with $\phi_\b \notin R_\b$
and thus it follows form the construction that $y_\b \phi_\b
\notin G_{\b +1}$. On the other hand it follows from Lemma
\ref{desc2} that $y_\b \phi_\b = y_\b \phi \in G \cap
\widehat{P_\b} \subseteq G_{\b +1}$ - a contradiction. Thus $\End
G = R$. \fine

We conclude this section with an immediate corollary of Theorem
\ref{super}. In Section \ref{warmup} we noted that $R = \Z F$ for
some free group $F$ has only trivial zero devisors. Therefore
Theorem \ref{super} has the following

\begin{corollary} For any cardinal $\l = \mu^+$ with $\mu^{\aln} =
\mu$ there is an $\ale$-free, indecomposable abelian UT-group of
cardinality $\l$.
\end{corollary}

\section{Discussions and applications}\label{Discus}

We want to discuss some consequences of our Main Theorem
\ref{super}. Let $\Mon G$ be the monoid of all monomorphisms of an
abelian group $G$. Also let be $$\Mont G = \{\v\in \Mon G, G/\Im
\v \mbox{ is torsion }\}.$$

We immediately note a

\begin{proposition}\label{mont} If $G$ is an abelian group of type $0$, then
$\Mont G$ is a submonoid of $\Mon G$.\end{proposition}

\Proof If $\v,\psi\in \Mont G$, then $\v\psi\in\Mon G$. To see
that $G/\Im(\v\psi)$ is torsion, consider any $x\in G$. There are
$y\in G, 0\ne n\in \o$ with $nx =y\psi$; similarly $my =z\psi$,
hence $mnx = m(y\psi) = z(\v\psi)$ and $\v\psi\in \Mont G$. \fine

We notice at many steps of the proof of Main Theorem \ref{super}
that UT (with respect to $\Aut G$) is a very strong property. If
we replace $\Aut G$ by $\Mon G$ or $\End G$ in the definition of
UT, this is even stronger: We either get nothing new or arrive on
classical territory as shown next. We begin with an easy
\begin{observation}\label{monUT} If $G$ is an abelian group of type $0$ and
$\Mon G$ acts UT on $\pp G$, then $\Aut G =\Mon G$.
\end{observation}

\Proof Suppose $\v\in \Mon G$ and $x\in \pp G$, then $y=^\df
x\v\in \pp G$ and there is $\psi\in \Mon G$ with $y\psi = x$,
hence $x\v\psi =x, y\psi\v = y$ and $x\id =x, y\id = y$. Thus
$\v\psi = \psi \v = \id$ follows by the U-property and $\v\in \Aut
G$. \fine

Hence UT is the same for automorphisms and monomorphism. If we
require even more, that $\End G$ acts transitive on $G$, then $G$
is a vector space and the problem on the existence of UT-modules
becomes trivial.

\begin{definition} We will say that $G$ is E-transitive if for
any pair $x,y\in \pp G$ there is $\s\in\End G$ such that $x\s =
y$.
\end{definition}

Any $\s\in \Aut G$ induces a permutation of $\pp G$, but this does
not hold for $\s\in\End G$. Our next results removes the set
theoretic assumption V=L from a main theorem in Dugas, Shelah
\cite{DS} and strengthens the outcome. Recall from the
introduction that $A$-transitive in the sense of \cite{DS} is the
same as transitive (or just T) in this context.

\begin{corollary} For any cardinal $\l = \mu^+$ with $\mu^{\aln} =
\mu$ there is an $\ale$-free abelian group of cardinality $\l$
which is E-transitive, but not transitive.
\end{corollary}

In view of Proposition \ref{mont} we will strengthen this
corollary further and sketch a proof:
\begin{theorem} For any cardinal $\l = \mu^+$ with $\mu^{\aln} =
\mu$ there is an $\ale$-free abelian group $G$ of cardinality $\l$
with $\Aut G =\pm 1$, which is transitive with respect to $\Mont
G$. From $\Aut G = \pm 1$ follows immediately that $G$ is not a
T-group.
\end{theorem}

{\sl Sketch of a proof.} The proof is very similar to the proof
given for $\Aut G$ but simpler because the U-property must go, see
Observation \ref{monUT}. We must run through the paper once more,
essentially replacing $\Aut G$ by $\Mont G$. First we will replace
the definition $\pAut G$ by
$$\pmont G =^\df \{\v , \v:\Dom\v \arr \Im\v \subseteq G,\ \v
\mbox{ an isomorphism, } G/\Dom\v \ \ale\mbox{-free }\}.$$ If
$\FF\subseteq \pmont G$, then $\FFF$ is the submonoid of $\pmont
G$ which is generated as a monoid by $\{\pm 1\}\cup \FF$, hence
$\FFF\subseteq \pmont G$ by Proposition \ref{mont}. The Definition
\ref{partfree} of an $\ale$-free pairs $(G,(\v_t)_{t\in u})$
remains the same, except that $\pAut$ and automorphisms must be
replaced by $\pmont$ and monomorphisms, respectively. Moreover
condition $(iii)$ is weaker (automatically) because $\FFF$ is not
closed under weak inverses.

The U-property must be replaced by a very weak condition assuring
that finally a free monoid will act on the group as monomorphisms
which are not automorphisms. This is incorporated in the new
definition for $\KK$ (compare Definition \ref{maindef}). We say
that a quadruple $\x = (G,\FF,\pi,h)$ belongs to $\KK$ if and only
if the following holds.
\begin{enumerate}
\item  $G$ is an $\ale$-free abelian group.
\item $\FF = \{\v_t : t\in J\}$ is a family of symbols $\v_t$ indexed
by $J = J^\x$, which generates a free monoid $\FFF$.
\item $\pi:\FF \arr \pmont G$ is a map which naturally extends
to $\FFF$ and satisfies the weak U-property: If $\v, \v'\in \FFF$,
$D =\Dom\pi(\v)\cap \Dom \pi(\v') \ne 0$ and $x\v= x\v'$ for all
$x\in \pp G\cap D$ then $\v =\v'$.
\item $h:\PP_{\aln}(J) \arr \Im(h)$ is a partial function from
$\Dom h \subseteq \PP_{\aln}(J)$ such that for $u\in \Dom
h\subseteq G$ the following holds.
\begin{itemize}
\item[(a)] $h(u)$ is a countable subgroup of $G$ and
$(\Dom \v^\x_t \cap h(u))\v^\x_t \subseteq h(u)$ for all $t\in u$
\item[(b)] If $\Gbar = G/h(u)$ and $\vbar^\x_t$ denote the
monomorphism induced by $\v^\x_t$, then $(\Gbar,(\vbar^\x_t)_{t\in
u})$ is $\ale$-free.
\end{itemize}
\end{enumerate}

The definition of an order on $\KK$ remains the same, except that
condition $(iii)(d)$ must be removed. Following the arguments
along Section \ref{partialauto} we note that Lemma~ \ref{morepa}
and a simple version of Lemma \ref{ext} are crucial for the weak
U-property. We obtain a theorem parallel to Theorem \ref{Auto}:
For each $\x = (G,\FF,\pi,h) \in\KK$ we can find $\x' =
(G',\FF',\pi',h')\in\KK$ with $\x\le \x'$ such that $\FF'\subset
\pmont G'$ which freely generates $\FFFs$ as a submonoid, $\FFFs$
acts transitive on $\pp G'$, all members of $\FFF\setminus \{1\}$
are proper monomorphisms  and $|G'|+|\FF'|= |G|+ |\FF|$. As in
Section \ref{constrUT} we modify $G'$ and get a new $G$ with
endomorphism ring $R= \Z\FFF$. Recall that $\FFF$ is a free monoid
generated by $\FF$, hence $U(R) =\pm 1$. Modification of
Proposition \ref{notonto} will show that any $\v\in R\setminus
U(R)$ is not in $\Mon_t(G)$; thus $\Mon_t(G)= \FFF$ and $\Aut G =
U(R) = \pm 1$. The group $G$ cannot be transitive. \fine

\noindent
Address of the authors:\\
R\"udiger G\"obel \\ Fachbereich 6, Mathematik, Universit\"at
Essen, 45117 Essen, Germany \\ {\small e--mail:
R.Goebel@Uni-Essen.De}\\ and \\ Saharon Shelah \\ Department of
Mathematics, Hebrew University, Jerusalem, Israel
\\ and Rutgers University, Newbrunswick, NJ, U.S.A \\ {\small
e-mail: Shelah@math.huji.ac.il}
\end{document}